\documentclass{tran-l}
\usepackage{amssymb}
\usepackage{array}

 \newcounter{itemprimed} 
 \newcounter{itemdoubleprimed}

\newcommand{\tf}{\tilde{f}}
\newcommand{\M}[1]{\Gamma_{#1} \backslash G_{#1}}
\newcommand{\Mp}[1]{\Gamma'_{#1} \backslash G_{#1}}
\newcommand{\addtorus}[1]{#1^{\mathsf{T}}}

\newcommand{\Lie}{\mathcal}
\newcommand{\Ad}{\operatorname{Ad}}
\newcommand{\SO}{\operatorname{SO}}
\newcommand{\GL}{\operatorname{GL}}
\newcommand{\Rad}{\operatorname{Rad}}

\newcommand{\iso}{\cong}
\newcommand{\homeo}{\approx}

\newcommand{\prodsemi}{\rtimes}
\newcommand{\semiprod}{\ltimes}

\newcommand{\real}{\mathord{\mathbb R}}

 \newcommand{\Gs}{\addtorus{G_2}}

\newcommand{\quot}{\bar}
\newcommand{\closure}{\overline}
\newcommand{\G}{\closure{\Ad G}}
\newcommand{\Gi}[1]{\closure{\Ad G_{#1}}}
\newcommand{\Z}[1]{\closure{\Ad_G #1}}
\newcommand{\Zi}[2]{\closure{\Ad_{G_{#1}} #2}}
\newcommand{\SOc}{\operatorname{\cover{\SO}(2)}}
\newcommand{\cover}{\widetilde}

\newcommand{\fol}{\mathcal F}

\renewcommand{\ell}[1]{\tau_{#1}}  

\newcommand{\stackarrow}[1]{\setbox0\hbox{ \ $\scriptstyle#1$ \ }
 \setbox1\hbox{$\rightarrow$}\setbox2\hbox to \wd0{\rightarrowfill}\ht2=\ht1
 \mathrel{\mathop{\copy2}\limits^{\copy0}}}

\newcommand{\eqcomment}[1]{&& \qquad \text{#1}}

\renewcommand{\see}[1]{\textnormal{(}see~\textnormal{\ref{#1})}}

\renewcommand{\eqref}[1]{Eq.~\textnormal{(\ref{#1})}}

\swapnumbers

\numberwithin{equation}{section}

\theoremstyle{plain}

\newtheorem{thm}[equation]{Theorem}
\newtheorem{mainthm}[equation]{Main Theorem}
\newtheorem{lem}[equation]{Lemma}
\newtheorem{cor}[equation]{Corollary}
\newtheorem{prop}[equation]{Proposition}

\theoremstyle{definition}

\newtheorem{defn}[equation]{Definition}

\theoremstyle{remark}

\newtheorem{rem}[equation]{Remark}
\newtheorem{ack}[equation]{Acknowledgment}  
\newtheorem{eg}[equation]{Example}
\newtheorem{warn}[equation]{Warning}

 \newcounter{step}
 \newenvironment{step}[1][\unskip]{\refstepcounter{step}\em
 \medskip \noindent Step \thestep\ #1.}{\unskip\upshape}
 \renewcommand{\thestep}{\arabic{step}}

\begin{document}


\renewcommand{\datename}{Submitted to \emph{Trans.\ Amer.\ Math.\ Soc.}}

\makeatletter
\def\@serieslogo{\vbox to\headheight{%
        \vss}}
\def\@setcopyright{}

\makeatother

\mbox{ }\vskip-0.5in 


\title{Foliation-preserving maps between solvmanifolds}

\author{Holly Bernstein}
 \address{Department of Mathematics, Williams College, Williamstown, MA 01267}
 \curraddr{Department of Mathematics, Washington University, St.~Louis, MO 63130}
 \email{holly@math.wustl.edu}

\author{Dave Witte}
 \address{Department of Mathematics, Williams College, Williamstown, MA 01267}
 \curraddr{Department of Mathematics, Oklahoma State University, Stillwater, OK
74078}
 \email{dwitte@math.okstate.edu,
  http://www.math.okstate.edu/\char'176dwitte} 


\date{(February 1998).}

\begin{abstract}
 For $i = 1,2$, let $\Gamma_i$ be a lattice in a simply connected, solvable Lie
group~$G_i$, and let $X_i$ be a connected Lie subgroup of~$G_i$. The double cosets
$\Gamma_igX_i$ provide a foliation~$\fol_i$ of the homogeneous space
$\M{i}$. Let $f$ be a continuous map from $\M1$
to $\M2$ whose restriction to each leaf of~$\fol_1$ is a covering
map onto a leaf of~$\fol_2$. If we assume that $\fol_1$ has a dense leaf, and make 
certain technical technical assumptions on the lattices $\Gamma_1$ and~$\Gamma_2$,
then we show that $f$ must be a composition of maps of two basic types: a
homeomorphism of $\M1$ that takes each leaf of~$\fol_1$ to itself, and a map that
results from twisting an affine map by a homomorphism into a compact group.

 We also prove a similar result for many cases where $G_1$ and~$G_2$ are neither
solvable nor semisimple.
 \end{abstract}

\maketitle

\section{introduction}

Let $\Gamma_1$ be a lattice in a simply connected, solvable Lie group~$G_1$. Any
connected Lie subgroup~$X_1$ of~$G_1$ acts by translations on the homogeneous space
$\M1$; the orbits of this action are the leaves of a
foliation~$\fol_1$ of $\M1$. We call this \emph{the foliation of\/
$\M1$ by cosets of~$X_1$}. Now suppose $\Gamma_2$ is a lattice in some
other simply connected, solvable Lie group~$G_2$, and that $X_2$ is a connected Lie
subgroup of~$G_2$, with corresponding foliation~$\fol_2$ of $\M2$. It is natural to
ask whether $\fol_1$ is topologically equivalent to~$\fol_2$, or, more generally,
whether there is a continuous map~$f$ from $\M1$ to
$\M2$ whose restriction to each leaf of~$\fol_1$ is a covering map
onto a leaf of~$\fol_2$. If so, it is of interest to know all the possible maps~$f$.

Under the assumption that some leaf of~$\fol_1$ is dense, and technical assumptions
on the lattices $\Gamma_1$ and~$\Gamma_2$, we show that every possible~$f$ is a
composition of maps of the basic types described in Example~\ref{map-egs} below.
(Remarks \ref{reduce-to-conn} and~\ref{reduce-to-dense} show that there are always
finite covers of $\M1$ and $\M2$ that satisfy the technical assumptions
on the lattices.) The reader may note that the composition of maps of the types described
in \ref{map-egs}\ref{map-translate} and~\ref{map-egs}\ref{map-affine} is an
affine map; the composition of types \ref{map-egs}\ref{map-translate}
and~\ref{map-egs}\ref{map-almostdoubly} is a \emph{doubly crossed} affine map
(c.f.~\cite[Defn.~7.3]{Witte-super}).

\begin{eg} \label{map-egs}
 \mbox{ } 

 \begin{enumerate}

 \renewcommand{\theenumi}{\Alph{enumi}}

\item \label{map-leaftoself}
 If $\fol_1 = \fol_2$, let $f$ be a homeomorphism of $\M1$ that maps each leaf
of~$\fol_1$ to itself.

\item \label{map-translate}
 If $\M1 = \M2$, and $X_1 = r^{-1} X_2 r$ is
conjugate to~$X_2$, let $f$ be a translation: $f(\Gamma_1 g) = \Gamma_1 gr^{-1} $.

\item \label{map-affine}
 If there is a continuous group homomorphism $k \colon G_1 \to G_2$ such that
$k(\Gamma_1 ) \subset \Gamma_2$, and the restriction $k|_{X_1}$ of~$k$ to~$X_1$ is
an homeomorphism onto~$X_2$, let $f\colon \M1 \to \M2$ be the map induced by~$k$:
$f(\Gamma_1 g) = \Gamma_2 k(g)$.

 \item[\ref{map-affine}$'$.]
 \refstepcounter{itemprimed}
 \label{map-almostdoubly}
 A map  $f \colon \M1\to \M2$ of type~\ref{map-affine} can
usually be modified as follows. Embed $G_2$ as a closed subgroup of some solvable
Lie group~$\Gs$. For $i=1,2$, let $T_i$ be a compact, abelian subgroup of~$\Gs$,
and let $\delta_i \colon G_1 \to T_i$ be a homomorphism, such that $\delta_i(\Gamma_1
) = e$. Define $\phi \colon G_1 \to \Gs$ by 
 $ \phi(g) = k(g) \cdot \delta_1(g) \cdot \delta_2(g) $.
 Under appropriate hypotheses \see{phi(G)isgroup}, $\phi(X_1)$ is a subgroup of~$\Gs$
(even though $\phi$ is usually not a group homomorphism), and the restriction of~$\phi$
to each coset of~$X_1$ is a homeomorphism onto a coset of~$\phi(X_1)$.

Let $G_2'$ be any connected Lie subgroup of~$\Gs$ that contains $\phi(G_1)$, and
let $\Gamma_2'$ be a lattice in~$G_2'$. Then the cosets of the subgroup
$\phi(X_1)$ provide a foliation~$\fol_2'$ of $\Gamma_2' \backslash G_2'$. Assume
$\Gamma_2 \subset \Gamma_2'$, so $\phi$ induces a well-defined map $f_2 \colon
\M1 \to \Gamma_2' \backslash G_2'$ defined by $f_2(\Gamma_1 g) = \Gamma_2'
\phi(g)$. The restriction of~$f_2$ to each leaf of~$\fol_1$ is a covering map
onto a leaf of~$\fol_2'$. 

One could add more homomorphisms $\delta_3$, $\delta_4$, etc., but
Theorem~\ref{solv-fol-map} shows that this is not necessary.

 \item[\ref{map-affine}$''$.]
 \refstepcounter{itemdoubleprimed}
 \label{map-doublycrossed}
 Instead of assuming that $\delta_i(\Gamma) = e$, the construction described
in~\ref{map-almostdoubly} can still be carried out if we make the weaker
assumption that $\delta_1(\gamma)\delta_2(\gamma) = e$, for all $\gamma
\in\Gamma$.

\end{enumerate}

\end{eg}

 The precise statement of our result requires the definition of the almost-Zariski
closure of a subgroup.

\begin{defn}[{\cite[Defn.~3.2]{Witte-super}}]
 A subgroup~$A$ of $\GL_n(\real)$ is \emph{almost Zariski closed} if
there is a Zariski closed subgroup~$B$ of $\GL_n(\real )$, such that $B^\circ
\subset A \subset B$, where $B^\circ$ is the identity component of~$B$ in the
topology of $\GL_n(\real)$ as a $C^\infty$~manifold (not the Zariski
topology). There is little difference between being Zariski closed and almost
Zariski closed, because $B^\circ$ always has finite index in~$B$.
 \end{defn}

\begin{defn}[{\cite[Defn.~3.6]{Witte-super}}] \label{almost-def}
 Let $A$ be a subgroup of $\GL_n(\real)$. The \emph{almost-Zariski
closure}~$\closure{A}$ of~$A$
 is the unique smallest almost-Zariski closed subgroup that contains~$A$.
 In particular, if $A$ is a subgroup of a Lie group~$G$, we use $\Z{A}$ to
denote the almost-Zariski closure of $\Ad_G A$ in $\GL(\Lie G)$, where $\Lie G$
is the Lie algebra of~$G$.
 \end{defn}

\begin{mainthm} \label{solv-fol-map}
 Let $X_1$ and $X_2$ be connected Lie subgroups of simply connected, solvable Lie
groups $G_1$~and~$G_2$, respectively. 
 For $i=1,2$, let $\Gamma_i$ be a lattice in $G_i$. Assume that $\Zi{1}{\Gamma_1}
= \Gi{1}$, and that $\Zi{2}{\Gamma'}$ is connected, for every subgroup~$\Gamma'$
of~$\Gamma_2$. Assume, furthermore, that the foliation of~$\M1$ by cosets of~$X_1$
has a dense leaf.

 Let~$f$ be a continuous map from $\M1$ to $\M2$, such that the
restriction of~$f$ to each leaf of the foliation of~$\M1$ by cosets of~$X_1$ is a
covering map onto a leaf of the foliation of~$\M2$ by cosets of~$X_2$.
 Then there exists 
 \begin{itemize}
    \item a map $a \colon \M1 \to \M1$ of type~\ref{map-egs}\ref{map-leaftoself};
    \item a map $b \colon \M2 \to \M2$ of type~\ref{map-egs}\ref{map-translate};
and
    \item a map $c \colon \M1 \to \M2$ of
type~\ref{map-egs}\ref{map-doublycrossed},
  \end{itemize}
 such that $f(\Gamma _1 g) = b\Bigl( c \bigl( a(\Gamma g) \bigr) \Bigr)$,
 for all $g \in G_1$.

In the definition of~$c$, we may take $\Gs$ as defined in
Definition~\ref{ell-defn}. We may take $T_1$ to be the elliptic part of~$k(X_1)$
and let $T_2$ be an appropriately chosen elliptic part of~$r^{-1} X_2 r$, where
$k$ is the homomorphism used in the construction of~$c$, and $r$ is the element
of~$G_2$ used in the construction of~$b$.
 \end{mainthm}

\begin{defn}[c.f.~Proposition~\ref{ell-prop}] \label{ell-defn}
 Let $G$ be a simply connected, solvable Lie group, and let $T_G$ be a maximal
compact torus of $\G$. Define $\addtorus{G} = G \prodsemi T_G$.

For any connected Lie subgroup~$X$ of~$\addtorus{G}$, there is a compact, abelian
subgroup~$T_X$ of~$\addtorus{G}$, such that $\Ad_G T_X$ is a maximal compact torus of
$\Z{X}$. The subgroup~$T_X$ is the \emph{elliptic part} of~$X$; it is unique up to
conjugation by an element of~$X$. 

The \emph{nonelliptic part} of~$X$ is the unique simply connected Lie subgroup~$Y$
of~$\addtorus{G}$ such that $X T_X = Y T_X$ and $\Z{Y}$ has no nontrivial compact subgroup.
 \end{defn}

The theorem was proved by D.~Benardete~\cite[Thm.~A(b)]{Benardete} in the special
case where $X_1$ and~$X_2$ are one-dimensional, the map~$f$ is a homeomorphism, and
the almost-Zariski closures $\Gi1$ and $\Gi2$ have no nontrivial compact
subgroups. (However, he proved only that \emph{some} foliation-preserving
homeomorphism is a composition of the standard types, not that \emph{all} are.)
D.~Witte~\cite[Thm.~5.1]{Witte-fols} removed the dimension restriction on the
subgroups $X_1$ and~$X_2$, and replaced it with the weaker hypothesis that they
are unimodular. We use the same methods as Benardete and Witte. The map $\delta^*$
does not appear in the conclusions of \cite{Benardete} and~\cite{Witte-fols},
because $T_1$ and~$T_2$ must be trivial if $\Gi2$ has no compact subgroups. 

If the  foliation of $\M1$ is not assumed to have a dense leaf, then it is not
possible to obtain such a precise global conclusion about the form of~$f$. However,
the proof shows that there is a homomorphism $k \colon G_1 \to \Gs$ with
$k(\Gamma_1) \subset \Gamma _2$, such that $k(X_1)$ and $r^{-1} X_2 r$ have the
same nonelliptic part, for some $r \in G_2$.

D.~Benardete and S.~G.~Dani \cite{BenDani} recently provided families of examples
$G$, $\Gamma $, and~$X_1$, such that $\Z{\Gamma} \not = \G$, yet, if the foliation
of $\M1$ by cosets of~$X_1$ is topologically equivalent to the foliation of $\M1$
by cosets of~$X_2$, then $X_1$ is conjugate to~$X_2$. The foliations are
topologically equivalent to linear foliations of ordinary tori (by applying
Remark~\ref{reduce-to-dense}), but not via affine maps.

The previous work of D.~Benardete~\cite{Benardete} and D.~Witte~\cite{Witte-fols}
requires $G_1$ and~$G_2$ to be either solvable or semisimple. This is because the
proofs rely on the Mostow Rigidity Theorem, which, until recently, was only
known in these cases. Now that results of this type have been generalized to other
groups \cite[\S9]{Witte-super}, the proof can be generalized.  Therefore, in
the final section of this paper, we sketch an application of Benardete's method to
many groups that are neither solvable nor semisimple.  However, unlike our work
in the solvable case, our results in this general setting are not at all definitive,
because we impose severe restrictions on the subgroups~$X_1$ and~$X_2$. (However,
the restrictions are automatically satisfied if $X_1$ and~$X_2$ are one-dimensional.) We
have not attempted to push these methods to their limit, because it seems clear that new
ideas will be needed to settle the general case.

\begin{ack}
 This paper grew out of the first author's undergraduate honors thesis (Williams
College, 1993).

 The second author would like to acknowledge that ideas of the late Louis
Auslander are the foundation of all of his research on homogeneous spaces,
including this paper. He was partially supported by research grants from the
National Science Foundation.
 \end{ack}

\section{Preliminaries}

\subsection{Technical assumptions on the lattices}
 The following remarks show that the assumptions on the lattices $\Gamma _1$
and~$\Gamma _2$ in the statement of Theorem~\ref{solv-fol-map} can be satisfied by
passing to finite covers of $\M1$ and $\M2$. Therefore, modulo finite covers, the
theorem describes the foliation-preserving maps for the natural foliations of all
solvmanifolds.

\begin{rem} \label{reduce-to-conn}
  The assumption in Theorem~\ref{solv-fol-map} that $\Zi{2}{\Gamma'}$ is
connected, for every subgroup~$\Gamma'$ of~$\Gamma_2$, can always be satisfied by
replacing $\Gamma_2$ with a finite-index subgroup (c.f.~\cite[Thm.~6.11,
p.~93]{Raghunathan}), or, in other words, by passing to a finite cover of $\M2$.
(This may also require $\Gamma _1$ to be replaced by a finite-index subgroup, so
that the map~$f$ is still well-defined.) However, the proof of the theorem does
not require the full strength of even this mild assumption. Rather, there is one
particular subgroup~$\Gamma'$ whose almost-Zariski closure needs to be connected;
see the first paragraph of the proof of the theorem. In particular, if $f$ is a
homeomorphism, then we need only assume $\Zi{2}{\Gamma_2}$ is connected.
 \end{rem}

\begin{rem} \label{reduce-to-dense}
 The assumption in Theorem~\ref{solv-fol-map} that $\Zi1{\Gamma_1} = \Gi1$ is
restrictive, but it does not limit the applicability of the result very severely,
because the theorem applies to a certain natural finite cover $\Gamma\backslash G$
of~$\M1$, which we now describe. (Note, however, that the covering map is usually
not affine; $G$ is \emph{not} isomorphic to~$G_1$.) Because $X_1$ has a dense
orbit on $\M1$, it is easy to see that $\Zi1{X_1}$ contains a compact torus~$T$ of
$\Gi1$, such that $\Zi1{\Gamma _1} T = \Gi1$. Therefore, the nilshadow construction
(cf.~\cite[Prop.~8.2]{Witte-super}) yields a simply connected, normal subgroup~$G$
of $G_1 \prodsemi T$, such that 
 \begin{itemize}
    \item $G$ contains a finite-index subgroup~$\Gamma $ of~$\Gamma _1$, 
    \item $\Z{\Gamma } = \G$, and
    \item $G T$ = $G_1 T$.
 \end{itemize}
 Define $\Delta \colon G_1 \to G$ by $\Delta(g) \in gT$. Then $\Delta$ is a
homeomorphism, and $\Delta(\gamma g) = \gamma \Delta(g)$, for all $\gamma \in
\Gamma $ and $g \in G$. Therefore, $\Delta^{-1}$ induces a finite-to-one covering
map $\Delta^* \colon \M{} \to \M1$. Furthermore, because $T$ normalizes~$X_1$, we
see that, letting $X = \Delta(X_1)$, we have $\Delta(gX_1) = \Delta(g) X$, so
$X$ is a subgroup of~$G$, and $\Delta^*$ maps each leaf of the foliation of $\M{}$
by cosets of~$X$ to a leaf of the foliation of $\M1$ by cosets of~$X_1$. 
 \end{rem}

\subsection{Hypotheses needed for examples of type \ref{map-egs}\ref{map-doublycrossed}}
 The following lemma describes some simple hypotheses that guarantee the conditions
needed for the construction of examples of type~\ref{map-egs}\ref{map-almostdoubly}
or~\ref{map-egs}\ref{map-doublycrossed}.

\begin{lem} \label{phi(G)isgroup}
 For $i = 1,2$, let $G_i$ be a simply connected, solvable Lie group, let $T_i$ be a
compact, abelian subgroup of $\Gs$, and let $\delta_i \colon G_1 \to T_i$ be a
continuous group homomorphism.
 Let $k \colon G_1 \to \Gs$ be a continuous group homomorphism, and let $X_1$ be a
connected closed subgroup of~$G_1$, such that $k|_{X_1}$ is a homeomorphism onto
$k(X_1)$. 
 Define $\phi_i \colon G_1 \to \Gs$ by
 $$\phi_1(g) = k(g) \cdot \delta_1(g)
 \qquad \text{and} \qquad
 \phi_2(g) = k(g) \cdot \delta_1(g) \cdot \delta_2(g) .$$
 If 
 $$[k(G_1), T_1] \subset k( \ker \delta_1 \cap \ker \delta_2)$$
 and 
 $$[k(X_1), T_1] \subset k(X_1 \cap \ker \delta_1 \cap \ker \delta_2) ,$$
  then $\phi_1(G_1)$ and $\phi_1(X_1)$ are subgroups of~$\Gs$, and the restriction
of~$\phi_1$ to each left coset of~$X_1$ is a homeomorphism onto a left coset of
$\phi_1(X_1)$. If, furthermore, 
 $$ \text{$T_2$ normalizes $\phi_1(X_1)$} ,$$
 then $\phi_2(X_1)$ is a subgroup, and the restriction of~$\phi_2$ to each left coset
of~$X_1$ is a homeomorphism onto a left coset of~$\phi_2(X_1)$.
 \end{lem}

\begin{proof}
 We give here only the last part of the proof, showing that the restriction of~$\phi_2$
to each coset of~$X_1$ is a homeomorphism onto a coset of~$\phi_2(X_1)$, because the
rest is very similar. Given $g \in G_1$ and $x \in X_1$, we have
 \begin{alignat*}{2}
 \phi_2(gx)
 &= k(gx) \delta_1(gx) \delta_2(gx) \\
 &= k(g) k(x) \delta_1(g) \delta_1(x) \delta_2(g) \delta_2(x) 
   \eqcomment{($k,\delta_1,\delta_2$ are homomorphisms)} \\
 &= k(g) \delta_1(g) k(xy) \delta_1(x) \delta_2(g) \delta_2(x) 
   \eqcomment{($\exists y \in X_1 \cap \ker \delta_1 \cap \ker \delta_2$)} \\
 &= k(g) \delta_1(g) k(xy) \delta_1(xy) \delta_2(xy) \delta_2(g) 
   \eqcomment{($y \in \ker \delta_i$
       and $T_2$ is abelian)} \\
 &= k(g) \delta_1(g) \phi_2(xy) \delta_2(g) \\
 &= k(g) \delta_1(g) \delta_2(g) \phi_2(x') 
   \eqcomment{($T_2 \subset N\bigl(\phi_1(X_1) \bigr) \Rightarrow \exists x' \in X_1$)}
\\
 &= \phi_2(g) \phi_2(x') .
  \end{alignat*}
 So $\phi_2(gX_1) = \phi_2(g) \phi_2(X_1)$.

We now show that $\phi_2|_{gX_1}$ is one-to-one. Suppose $\phi_2(gx_1) = \phi_2(gx_2)$,
for some $g \in G_1$ and $x_1,x_2 \in X_1$. Let $x_1',x_2' \in X_1$ with
 $k(x_j)
\delta_1(g) = \delta_1(g) k(x_j')$ and 
 $\delta_i(x_j') = \delta_i(x_j)$ for $i,j = 1,2$.
 Then
 \begin{alignat*}{2}
 \phi_1(g) \phi_1(x_1') \delta_2(x_1') \delta_2(g)
 &= k(g) \delta_1(g) k(x_1') \delta_1(x_1') \delta_2(x_1') \delta_2(g) \\
 &= k(g) k(x_1) \delta_1(g) \delta_1(x_1) \delta_2(x_1) \delta_2(g) \\
 &= k(g x_1) \delta_1(g x_1) \delta_2(g x_1) \\
 &= \phi_2(g x_1) \\
 &= \phi_2(g x_2) \\
 &=  \phi_1(g) \phi_1(x_2') \delta_2(x_2') \delta_2(g)
   \eqcomment{(reverse steps).} 
 \end{alignat*}
 Therefore $\phi_1(x_1') \delta_2(x_1') = \phi_1(x_2') \delta_2(x_2')$. Because
$\phi_1|_{X_1}$ is a homeomorphism onto $\phi_1(X_1)$, we know that $\phi_1(X_1)$ is
simply connected, so it has no nontrivial compact subgroups \see{sc-no-cpct}. Thus,
$\phi_1(X_1) \cap T_2 = e$, so because $\phi_1(x_1') \delta_2(x_1') = \phi_1(x_2')
\delta_2(x_2')$, we conclude that $\phi_1(x_1') = \phi_1(x_2')$, so $x_1' = x_2'$. In
particular, we have $\delta_i(x_1) = \delta_i(x_1') = \delta_i(x_2') = \delta_i(x_2)$ for
$i = 1,2$, so the equation $\phi_2(gx_1) = \phi_2(gx_2)$ immediately implies $k(x_1) =
k(x_2)$. Therefore $x_1 = x_2$.

The restriction of~$k$ to~$X_1$ is proper  (since it is a homeomorphism onto
its image), so the restriction of~$\phi_2$ to $gX_1$ is proper. From the preceding
paragraph, we know that it is also injective. Therefore, it is a homeomorphism onto its
image \see{proper-homeo}.
 \end{proof}

\subsection{Elliptic and nonelliptic parts of a subgroup}

The following proposition justifies the assertions in Definition~\ref{ell-defn},
and establishes some basic facts that will often be used without specific
reference.

\begin{prop} \label{ell-prop}
 Let $G$ be a simply connected, solvable Lie group, let $T_G$ be a maximal compact
torus of $\G$, define $\addtorus{G} = G \prodsemi T_G$, and let $X$ be a connected Lie
subgroup of~$\addtorus{G}$. Then 
 \begin{enumerate}
    \item $[\addtorus{G},\addtorus{G}] = [G,G]$, 
    \item $Z(\addtorus{G})$ has no nontrivial, compact subgroups,
    \item there is a compact, abelian subgroup~$T_X$ of~$\addtorus{G}$, unique up to
conjugation by an element of~$X$, such that $\Ad_G T_X$ is a maximal compact torus
of $\Z{X}$,
    \item there is a unique closed, simply connected subgroup~$Y$ of~$\addtorus{G}$ such
that $X T_X = Y T_X$ and $\Z{Y}$ has no nontrivial compact subgroup,
    \item we have $[X T_X, X T_X]  \subset X \cap Y$,
    \item $T_X \cap Y = e$, 
   and 
    \item $Y$ is normal in~$X T_X$.
 \end{enumerate}
 \end{prop}

\begin{proof} Lemma~\ref{comm-semi} asserts that $[\addtorus{G},\addtorus{G}] = [G,G]$.

Every compact subgroup of~$\addtorus{G}$ is conjugate to a subgroup of~$T_G$
\see{cpct-conj}. Since $T_G$ is a subgroup of $\operatorname{Aut} G$, we know that
no nontrivial element of~$T_G$ centralizes~$G$. Therefore $Z(\addtorus{G})$ has no
nontrivial, compact subgroups.

All maximal compact tori of $\G$ are conjugate under $\Ad G$ (e.g., see
\cite[Cor.~4.22]{Witte-SolvTess}), so, replacing $T_G$ by a conjugate, we may
assume $T_G$ contains a maximal compact torus~$S$ of $\Z{X}$. Then the desired
subgroup~$T_X$ is simply~$S$, thought of as a subgroup of $T_G \subset \addtorus{G}$. The
uniqueness follows from the fact that all maximal compact tori of $\Z{A}$ are
conjugate under $\Ad A$ (e.g., see \cite[Cor.~4.22]{Witte-SolvTess}).

Assume, as in the preceding paragraph, that $T_G$ contains~$T_X$. There is a
natural projection from $\G$ to~$T_G$, given by the splitting $\G = (A
\times T_G) \semiprod U$, where $A$ is a maximal $\real $-split torus and
$U$ is the unipotent radical.
 Let $\sigma  \colon \addtorus{G} \to T_G$ be the composite homomorphism
 $$\sigma  \colon \addtorus{G}
 \stackarrow{\Ad} \G
 \stackarrow{\operatorname{projection}} T_G
  \stackarrow{x \mapsto x ^{-1}} T_G \; ,$$ 
 and define $\Delta  \colon \addtorus{G} \to \addtorus{G}$ by
$\Delta (g) = g \sigma(g)$, so $\Delta$ is a \emph{nilshadow map}
\cite[Defn.~4.1]{Witte-super}. Then $\Delta (\addtorus{G})$ is a subgroup of $\Delta
(\addtorus{G})$ \cite[Cor.~4.8]{Witte-super}. Since $\Delta $ is obviously a proper map, we
know that $\Delta (G)$ is closed. Since $T_X \subset T_G$, we have $\sigma (X)
\subset T_X$, so $\Delta (X)$ is a subgroup of~$\Delta (G)$
\cite[Cor.~4.9]{Witte-super} and, obviously, $AT_X = \Delta (X) T_X$. By
construction, $\Z{\Delta (\addtorus{G})}$ has no nontrivial compact subgroup
\cite[Prop.~4.10]{Witte-super}. Since $Z(\addtorus{G})$ has no nontrivial, compact
subgroups, this implies that $\Delta (\addtorus{G})$ has no nontrivial, compact subgroups,
which means that $\Delta (\addtorus{G})$ is simply connected \see{sc-no-cpct}. Therefore,
we may let $Y = \Delta (X)$ \see{subgrp-sc}. If $Y'$ is any nonelliptic part
of~$X$, then, because $\Z{Y'}$ has no nontrivial compact subgroup, the
subgroup~$Y'$ must be contained in the kernel of~$\sigma$. This kernel is
precisely $\Delta (\addtorus{G})$, so it is not difficult to see that $Y' = Y$.

The definition of the nilshadow map~$\Delta$ immediately implies $X \cap
\ker\sigma  \subset \Delta (X)$. Since $T_G$ is abelian, $\ker\sigma $ must
contain $[\addtorus{G},\addtorus{G}]$. Therefore $X \cap [\addtorus{G},\addtorus{G}] \subset Y$. Then, because
Lemma~\ref{comm-semi} implies $[XT_X,XT_X] = [X,X] \subset X$, we have
 $[XT_X,XT_X] \subset X \cap [\addtorus{G},\addtorus{G}] \subset X \cap Y$.

Being simply connected, $Y$ has no nontrivial, compact subgroups
\see{sc-no-cpct}, so $T_X \cap Y = e$.

Since $[Y , X T_X] \subset [X T_X , X T_X] \subset X \cap Y \subset Y$, we see that $Y$
is normal in $X T_X$.
 \end{proof}

\subsection{Nondivergent subgroups}

\begin{defn}[{\cite[\S4]{Witte-fols}}] \label{diverge-defn}
 Let $X$~and~$Y$ be subsets of a Lie group~$G$. We say that \emph{$X$ does not
diverge from~$Y$} if there is a compact subset~$K$ of~$G$ with $X \subset Y K$. 
 If $X$ does not diverge from~$Y$, and $Y$ does not diverge from~$X$, then we may say
that \emph{$X$ and~$Y$ do not diverge from each other.}
 \end{defn}

For the special case where the subgroup~$Y$ is unimodular, the following
proposition was proved by D.~Witte \cite[Cors.~4.10 and 4.11]{Witte-fols}.

\begin{prop} \label{no-cpct-diverge}
 Let $X$~and~$Y$ be connected Lie subgroups of a simply connected solvable Lie
group~$G$, and assume that $\G$ has no nontrivial compact subgroups. If $X$ does
not diverge from~$Y$, then $X \subset Y$. 
 \end{prop}

\begin{proof}
 Each element of $[G,G]$ acts unipotently, hence unimodularly, on each subspace
of~$\Lie G$ that it normalizes, so $[G,G] \cap Y \subset \ker \Delta_Y$, where
$\Delta_Y$ is the modular function of~$Y$. Hence $\Delta_Y$ extends to a
continuous homomorphism $\Delta  \colon G \to \real^+$. Define a semidirect
product $G \semiprod \real$, by letting each $g \in G$ act on~$\real$ via
multiplication by $1/\Delta (g)$, so $Y \semiprod \real$ is a unimodular subgroup
of $G \semiprod \real$. Since $X$ does not diverge from~$Y$ in~$G$, we see that
$X$ does not diverge from $Y \semiprod \real$ in $G \semiprod \real$, so
the proof of \cite[Cor.~4.11]{Witte-fols} shows that $X \subset Y \semiprod
\real$. (Although the statement of the corollary assumes both $X$ and~$Y$ are
unimodular, the proof only requires this assumption for~$Y$.) Hence $X  \subset (Y
\semiprod \real) \cap G = Y$, as desired.
 \end{proof}

\begin{cor} \label{diverge-nonell}
 Let $G$ be a simply connected, solvable Lie group, and let $X$~and~$Y$ be
connected Lie subgroups of~$\addtorus{G}$. If $X$ does not diverge from~$Y$, then the
nonelliptic part of~$X$ is contained in the nonelliptic part of~$Y$.
 \end{cor}

\begin{proof}
 Letting $X_1$, $Y_1$ and~$G_1$ be the nonelliptic parts of $X$, $Y$, and~$G$,
respectively, it is easy to see that $X_1$ does not diverge from~$Y_1$ in~$G_1$.
Therefore, Proposition~\ref{no-cpct-diverge} implies $X_1 \subset Y_1$, as desired.
 \end{proof}

\begin{eg} Let $G = \SOc \semiprod \real^2$, where $\SOc$ is the universal cover
of $\SO(2)$, and let $X$ be a subgroup of~$G$ that is conjugate to $\SOc$. Then $X$
does not diverge from $\SOc$ in~$\addtorus{G} \iso \real \times (\SO(2) \prodsemi
\real^2)$, but $X$ need not be contained in $\SOc$, so we see that the conclusion
of the corollary cannot be strengthened to say that $X \subset Y$.
 \end{eg}

\subsection{Miscellaneous facts}

For ease of reference, we record some basic results on solvable groups.

\begin{lem}[{\cite[Thm.~XII.2.2, p.~137]{Hoch}}] \label{subgrp-sc}
  Every connected subgroup of a simply connected, solvable Lie group~$G$ is closed
and simply connected.
 \end{lem}

\begin{lem}[{\cite[Lem.~2.17]{Witte-super}}] \label{sc-quot}
 Let $N$ be a closed subgroup of a connected, solvable Lie
group~$G$. Then $G/N$ is simply connected if and only if $N$ is connected and
contains a maximal compact subgroup of~$G$.
 \end{lem}

\begin{lem}[{\cite[Cor.~2.18]{Witte-super}}] \label{sc-no-cpct}
 A connected, solvable Lie group is simply connected if and only if it has no
nontrivial compact subgroups.
 \end{lem}

\begin{lem}[{cf.~\cite[Lem.~3.24]{Witte-super}}] \label{comm-semi}
 If $G$ is a connected, solvable Lie group, and $T$ is an abelian subgroup of
$\G$, then $[G \prodsemi T, G \prodsemi T] = [G,G] \prodsemi e$.
 \end{lem}

\begin{prop}[{\cite[Thm.~XV.3.1, pp.~180--181]{Hoch}}] \label{cpct-conj}
 Every compact subgroup of a connected Lie group~$G$ is contained in a maximal
compact subgroup, and all maximal compact subgroups of~$G$ are conjugate.
\end{prop}

\begin{lem} \label{torus-quot}
 Let $\Gamma_1$ be a lattice in a simply connected, solvable Lie group $G_1$, such
that $\Zi{1}{\Gamma_1} = \Gi{1}$, and let $E$ be a compact, abelian Lie group. 
 Let $X_1$ be a connected Lie subgroup of $G_1$, and assume that the foliation
of~$\M1$ by cosets of~$X_1$ has a dense leaf. Suppose $\ell{} \colon X_1 \to E$
and $\delta^* \colon \M1 \to E$ are continuous maps, such that $\delta^*(\Gamma_1)
= e$, and 
 $\delta^*(px) = \delta^*(p)\ell{}(x)$, for every $p \in \M1$ and $x \in X_1$.
 Then there is a continuous homomorphism $\cover{\delta} \colon G_1 \to E$ such
that 
 $\delta^*(\Gamma_1 g) = \cover{\delta}(g)$ for every $g \in G_1$.
 \end{lem}

\begin{proof}[Proof \upshape{(c.f.~\cite[pf.~on p.~502]{Benardete})}]
 For $x,y \in X_1$, we have 
 $$\delta^*(\Gamma_1) \ell{}(xy) = \delta^*(\Gamma_1 xy) 
 = \delta^*(\Gamma_1 x) \ell{}(y) = \delta^*(\Gamma_1) \ell{}(x) \ell{}(y) ,$$
 so we see that $\ell{}$ is a homomorphism. Because $X_1$ is simply connected
\see{subgrp-sc}, we may lift~$\ell{}$ to a homomorphism $\cover{\ell{}} \colon X_1
\to \cover{E}$, where $\cover{E}$ is the universal cover of~$E$.

  Because the foliation of $\Gamma_1 \backslash G_1$ has a dense leaf, we may
assume $\Gamma_1 X_1$ is dense in~$G_1$. Because $G_1$ is simply connected, we
may lift $\delta^*$ to a map $\cover{\delta} \colon G_1 \to \cover{E}$ with
$\cover{\delta}(e) = e$. The restriction of~$\cover{\delta}$ to the fundamental
group~$\Gamma_1$ is a homomorphism into~$\cover{E}$. Theorem~\ref{super} implies
that this restriction $\cover{\delta}|_{\Gamma_1}$ extends to a continuous homomorphism
$k \colon G_1 \to \cover{E}$.

 We have $\cover{\delta}(\gamma g) = k(\gamma) \cdot \cover{\delta}(g)$, for every
$\gamma \in \Gamma_1$ and $g \in G_1$. Therefore, because $\M1$ is compact, there
is a compact subset~$K$ of~$\cover{E}$, such that $\cover{\delta}(g) \in k(g) K$,
for all $g \in G_1$. In particular, for $x \in X_1$, we have $\cover{\ell{}}(x) =
\cover{\delta}(x) \in k(x) K$, so the difference $\cover{\ell{}} - k$ is a
homomorphism with bounded image. Therefore, Lemma~\ref{sc-no-cpct} implies that
the image is trivial, which means $\cover{\ell{}} = k$, so $\cover{\delta}$ agrees
with~$k$ on~$X_1$. Since they also agree on~$\Gamma_1$, and $\Gamma _1 X_1$ is
dense in~$G_1$, this implies that $\cover{\delta} = k$ is a homomorphism.
 \end{proof}

In the statement of the following result in \cite{Witte-super}, it is assumed that
the maximal compact torus~$T_{G_2}$ used in the construction of~$\Gs$ contains a
maximal compact torus of $\Zi2{\Gamma_1^\alpha }$. Because all maximal compact
tori of $\Gi2$ are conjugate under $\Ad G_2$, this assumption is unnecessary.

\begin{thm}[{\cite[Cor.~6.5]{Witte-super}}] \label{super}
 Let $\Gamma_1$ be a lattice in a simply connected, solvable Lie group~$G_1$,
and assume $\Zi1{\Gamma_1} = \Gi1$. 
 Let $G_2$ be a simply connected, solvable Lie group. If $\alpha$ is a
homomorphism from~$\Gamma_1$ into $G_2$, such that
 $\Zi2{\Gamma_1^\alpha }$ is connected,
 then $\alpha$ extends to a continuous homomorphism from~$G_1$ to~$\Gs$.
 \end{thm}

For convenience, we also note the following well-known, simple lemma.

\begin{lem} \label{proper-homeo}
 Every continuous, proper bijection between locally compact Hausdorff
topological spaces is a homeomorphism.
 \end{lem}

\section{Proof of the Main Theorem}

The outline of this proof is based on \cite[\S6]{Witte-fols}. However, 
complications are caused by the possible lack of an inverse to~$f$, and by the
possible existence of nontrivial compact subgroups of $\Gi2$.

\begin{proof}[Proof of Theorem~\ref{solv-fol-map}]
 By composing~$f$ with the translation by some element $r \in G_2$, we may assume
without loss of generality that $f(\Gamma_1) = \Gamma_2$. Then, because $G_1$ is
simply connected, we may lift $f$ to a map $\tf\colon G_1 \to G_2$ with $\tf(e) =
e$. Because $\Gamma_i$ is the fundamental group of $\M{i}$, we see that the
restriction of~$\tf$ to~$\Gamma_1$ is a homomorphism into~$\Gamma_2$. Because
$\Zi2{\tf(\Gamma_1)}$ is connected, Theorem~\ref{super} implies that this
restriction $\tf|_{\Gamma_1}$ extends to a continuous homomorphism $k \colon G_1 \to
\Gs$.

\begin{rem} \label{bdd-dist}
 We have $\tf(\gamma g) = k(\gamma) \cdot \tf(g)$, for every $\gamma \in \Gamma_1$ and
$g \in G_1$. Therefore, because $\M1$ is compact, there is a compact subset~$K$
of~$\Gs$, such that $\tf(g) \in k(g) K$, for all $g \in G_1$. Hence, for every
subset~$A$ of~$G$, the sets $\tf(A)$ and $k(A)$ do not diverge from each other (see
Definition~\ref{diverge-defn}).
 \end{rem}

\begin{step} \label{f-homeo}
 The restriction of~$\tf$ to each coset of~$X_1$ is a homeomorphism onto a
coset of~$X_2$.
 \end{step}
 By assumption, the restriction of~$f$ to each leaf of the foliation of~$\M1$ by cosets
of~$X_1$ is a covering map onto a leaf of the foliation of~$\M2$ by cosets of~$X_2$.
Therefore, the restriction of~$f$ to each coset of~$X_1$ is a covering map onto a coset
of~$X_2$. Because $X_1$~and~$X_2$ are simply connected \see{subgrp-sc}, this covering
map must be a homeomorphism.

\begin{step}[{\cite[Step~3 of pf.~of Thm.~6.1]{Witte-fols}}] \label{same-nonell}
 $k(X_{1})$ and $X_{2}$ have the same nonelliptic part; call it~$Y$.
 \end{step} 
 Because $k$ maps $X_1$ to~$k(X_1)$, and $\tf$ maps $X_1$ to~$X_2$,
Remark~\ref{bdd-dist} implies that $k(X_{1})$ and $X_{2}$ do not diverge from each
other. Therefore, Corollary~\ref{diverge-nonell} implies that $k(X_{1})$ and
$X_{2}$ have the same nonelliptic part.

\begin{step} \label{k-homeo}
 The restriction of~$k$ to~$X_1$ is a homeomorphism onto $k(X_1)$, and $k(X_1)$ is
closed.
 \end{step}
 Since $\tf|_{X_1}$ is a homeomorphism onto~$X_2$, it is a proper map. Therefore,
Remark~\ref{bdd-dist} implies that $k|_{X_1}$ is also a proper map. This implies
$k(X_1)$ is closed.  It also implies that the kernel of $k|_{X_1}$ is compact. Then
Lemma~\ref{sc-no-cpct} implies that the kernel is trivial, so $k|_{X_1}$ is
injective. Thus, $k|_{X_1}$ is an isomorphism onto its image \cite[Lem.~2.5.3,
p.~59]{Varadarajan}.

\begin{step}[{\cite[Step~4 of pf.~of Thm.~6.1]{Witte-fols}}]
\label{delta-in-N_G(Y)}
 For $g \in G_1$, define $\delta(g) \in \Gs$ by: $\tf(g) = k(g) \cdot
\delta(g)$; then $\delta(g)$ normalizes~$Y$, and $\delta(\gamma g) = \delta(g)$, for
every $\gamma \in \Gamma_1$, so $\delta$ factors through to a well-defined map $\delta^*
\colon \M1 \to N_{\Gs}(Y)$.
 \end{step}
 Because $\tf$ maps cosets of~$X_{1}$ to cosets of~$X_{2}$, we have
 $$ k(g)^{-1} \tf(gX_{1}) = k(g)^{-1} \tf(g) X_{2} = \delta(g) X_{2} = \bigl( \delta(g)
X_{2} \delta(g)^{-1} \bigr) \cdot \delta(g) . $$
 Then, because Remark~\ref{bdd-dist} implies that $k(g)^{-1} \tf(gX_1)$ and $k(X_1)
= k(g)^{-1} k(gX_1)$ do not diverge from each other, this implies that the subgroups
$\delta(g) X_{2} \delta(g)^{-1}$ and $k(X_{1})$ do not diverge from each other.
Therefore, Corollary~\ref{diverge-nonell} implies that the nonelliptic part of
$\delta(g) X_{2} \delta(g)^{-1}$ is the same as the nonelliptic part of~$k(X_1)$,
namely,~$Y$. On the other hand, because the nonelliptic part of~$X_2$ is~$Y$ (see
Step~\ref{same-nonell}), it is obvious that the nonelliptic part of~$\delta(g) X_{2}
\delta(g)^{-1}$ is~$\delta(g) Y \delta(g)^{-1}$. Therefore, $Y = \delta(g) Y
\delta(g)^{-1}$.

Because $\tf(\gamma g) = k(\gamma) \cdot \tf(g)$, it is easy to see that
$\delta(\gamma g) = \delta(g)$.

\begin{step}[{\cite[Step~5 of pf.~of Thm.~6.1]{Witte-fols}}] \label{T1-T2-Y}
 Let $T_{1}$ and $T_{2}$ be the elliptic parts of $k(X_{1})$ and~$X_{2}$,
respectively; then $\delta(g) \in T_{1} T_{2} Y$, for every $g \in G_1$.
 \end{step}
  For $x_{1} \in X_{1}$,
we have $\tf(gx_{1}) = \tf(g)\cdot x_{2}$ for some $x_{2} \in X_{2}$, so
 $$ k(g) k(x_1) \delta(g x_1) = k(gx_{1}) \cdot \delta(gx_{1}) = \tf(gx_{1}) =
\tf(g)x_{2} = k(g) \cdot \delta(g) \cdot x_{2}. $$
 Writing $k(x_{1}) = t_{1}y_{1}$ and $x_{2} = t_{2}y_{2}$ for some $t_{1} \in
T_{1}$, $t_{2} \in T_{2}$, and $y_{1},y_{2} \in Y$, we then have $t_{1}y_{1} \cdot
\delta(gx_{1}) = \delta(g) \cdot t_{2} y_{2}$. This implies that the map
$\overline\delta \colon \M1 \to T_{1} \backslash N_{\Gs}(Y)/ T_{2} Y$, induced
by~$\delta$, is constant on each leaf of the foliation of $\M1$ by cosets
of~$X_1$. Because this foliation has a dense leaf, this implies that
$\overline\delta$ is constant. Because $\tf(e) = e = k(e)$, we know that
$\delta(e) = e$, so this implies that $\delta(g)$ belongs to $T_{1} T_{2} Y$ for
every $g \in G_1$, as desired.

\begin{step} \label{nu-homeo}
 We may assume $(T_1T_2) \cap Y = e$; then the maps
 $$ \renewcommand{\arraystretch}{1.5}
 \begin{array}{rrrll}
   \nu \colon &Y \times T_1 T_2 &\to& Y T_1 T_2 \colon &(y,t) \mapsto y t
 ,  \cr
   \nu_1 \colon & X_1 \times T_1 T_2 &\to& k(X_1) T_1 T_2 \colon & (x,t)
\mapsto k(x) t
 , 
 \qquad \text{and}
\cr
   \nu_2 \colon &X_2 \times T_2 T_1 &\to& X_2 T_2 T_1 \colon & (x,t) \mapsto
x t
   \cr
 \end{array} $$
 are homeomorphisms.
 \end{step}
 Let $S = T_1 \cap (T_2 Y)$. Then there is some $g \in T_2 Y$, such that $S \subset
g^{-1} T_2 g$ \see{cpct-conj}. We have
 $$ T_1 \cap \bigl( (g^{-1} T_2 g) Y \bigr) = T_1 \cap (T_2 Y) = S  \subset g^{-1} T_2 g
,$$
 so, by replacing the choice $T_2$ of the elliptic part of~$X_2$ with the equally valid
choice $g^{-1} T_2 g$, we may assume $T_1 \cap (T_2 Y) \subset T_2$. We now show that
this implies  $(T_1T_2) \cap Y = e$: if $t_1 t_2 \in Y$, with $t_1 \in T_1$ and $t_2 \in
T_2$, then $t_1 \in T_1 \cap (T_2 Y)  \subset T_2$, so $t_1 t_2 \in T_2 \cap Y = e$, as
desired.

The maps $\nu$, $\nu_1$, and~$\nu_2$ are obviously continuous, surjective, and
proper. (For the properness of~$\nu_1$, recall that $k|_{X_1}$ is a homeomorphism
onto $k(X_1)$.) Thus, it suffices to show that they are injective
\see{proper-homeo}. 

Suppose $y' t_1' t_2' = y t_1 t_2$, for some $y',y \in Y$ and $t_i',t_i \in T_i$.
Then $t_1^{-1} t_1' t_2' t_2^{-1} \in (T_1 T_2) \cap Y = e$, so $t_1' t_2' =
t_1 t_2$ and, hence, $y' = y$, so $\nu$ is injective. 

Suppose  $k(x') t_1' t_2' = k(x) t_1 t_2$, for some $x',x \in X_1$ and $t_i',t_i
\in T_i$. We have $k(x') = y't'$ and $k(x) = yt$, for some $y',y \in Y$ and $t',t
\in T_1$. Then $y't' t_1' t_2' = yt t_1 t_2$, so, from the conclusion of the
preceding paragraph, we must have $y' = y$. Since $k(X_1) \cap T_1 = e$, this
implies $k(x') = k(x)$. Therefore $x = x'$, so $\nu_1$ is injective. A similar
argument applies to~$\nu_2$.

\begin{warn}
 The compact set $T_1 T_2$ need not be a subgroup of~$\Gs$, because
$T_1$~and~$T_2$ need not commute with each other.
 \end{warn}

\begin{step} \label{phi-defn}
 There is a left\/ $\Gamma_1$-equivariant homeomorphism~$\phi$ of~$G_1$,
such that  $\tf (g) \in k\bigl(\phi (g) \bigr) T_{1} T_{2}$, for every $g \in
G_1$, and $\phi$ takes each left coset of~$X_1$ onto itself.
 \end{step}
  Since $\delta(g) \in T_{1} T_{2} Y = k(X_1) T_{1} T_{2}$, there is a unique
element $\chi(g)$ of~$X_1$, such that $\delta(g) \in k\bigl( \chi(g) \bigr)  T_{1}
T_{2}$; namely, $\chi(g)$ is the first coordinate of $\nu_1^{-1} \bigl( \delta(g)
\bigr)$, so $\chi$ is a continuous function of~$g$.
 Define $\phi(g) = g \cdot \chi(g)$, so $\phi  \colon G_1 \to G_1$ is continuous,
and takes each left coset $gX_1$ of~$X_1$ into itself. Then 
 \begin{equation}
 \tf(g) = k(g) \delta(g) \in k(g) k\bigl( \chi(g) \bigr)  T_{1} T_{2}
 = k \bigl( g \cdot \chi(g) \bigr)  T_{1} T_{2}
 = k\bigl(\phi (g) \bigr) T_{1} T_{2} .
 \label{phi-eqn}
 \end{equation}
 So all that remains is to show that $\phi$ is left $\Gamma_1$-equivariant and has
a continuous inverse.

 Note that, because $\delta(\gamma g) = \delta(g)$, we must have $\chi(\gamma g) =
\chi(g)$, for all $g \in G_1$ and $\gamma \in \Gamma_1$. Therefore
 $$ \phi (\gamma g) = (\gamma g) \cdot \chi(\gamma g) = \gamma \bigl( g \cdot
\chi(g) \bigr) = \gamma \phi (g) ,$$
 which is exactly what it means to say that $\phi$ is left $\Gamma_1$-equivariant.

Define 
 $$\zeta \colon G_1 \times X_1 \times T_2 T_1 \to G_1
\times ( X_2 T_2 T_1 )
 \text{ by }
 \zeta(g,x,t) = 
 \bigl( g, \tf(g)^{-1} \tf(gx) t \bigr)
 . $$
 Assume for the moment that $\zeta$ is a homeomorphism; that is, $\zeta^{-1}$ is
continuous. For $h \in G_1$, let $\chi'(h)$ be the $X_1$-component of $\zeta^{-1}
\bigl(h, \delta(h)^{-1} \bigr)$, and define $\psi(h) = h \cdot \chi'(h) \in h X_1$, so
$\psi \colon G_1 \to G_1$ is continuous. From the definitions of $\chi'$ and~$\psi$, we
have 
 $$ \delta(h)^{-1} \in \tf(h)^{-1} \tf \bigl( h \chi'(h) \bigr) T_2 T_1
 = \bigl( k(h) \delta(h) \bigr)^{-1} \tf \bigl( \psi(h) \bigr) T_2 T_1 ,$$
 which means $\tf \bigl( \psi(h) \bigr) \in k(h) T_1 T_2$. However, from
\eqref{phi-eqn} (with $g = \psi(h)$), we also have
 $$ \tf \bigl( \psi(h) \bigr) \in k \Bigl( \phi \bigl( \psi(h) \bigr) \Bigr) T_1 T_2 .$$
 Since $\nu_1$ is injective, we conclude that $h = \phi \bigl( \psi(h) \bigr)$.
From \eqref{phi-eqn}, we have $k \bigl( \phi(g) \bigr) \in \tf(g) T_2 T_1$, so, because
$\tf$ maps $h X_1$ homeomorphically onto a coset of~$X_2$, and $\nu_2$ is injective, we
see that $\phi$ is injective. Therefore, $\psi$ is the inverse of~$\phi$.

Thus, all that remains is to show that $\zeta$ is a homeomorphism. Because $\tf$
maps each coset of~$X_1$ homeomorphically onto a coset of~$X_2$, and $\nu_2$ is
bijective, it is easy to see that $\zeta$ is bijective. Therefore, it suffices
to show that $\zeta$ is proper \see{proper-homeo}. Thus, let $\{g_n\} \subset
G_1$, $\{x_n\} \subset X_1$, and $\{t_n\} \subset T_2 T_1$, such that $\{
\zeta(g_n, x_n, t_n) \}$ converges. From the first coordinate of $\{
\zeta(g_n, x_n, t_n) \}$, we see that $\{g_n\}$ converges. Because $T_2 T_1$ is
compact, we may also assume, by passing to a subsequence, that $\{t_n\}$ converges.
Then, by considering the second coordinate of $\{ \zeta(g_n, x_n, t_n) \}$, we
see that $\{ \tf(g_n x_n) \}$ converges and, therefore, is bounded, so
Remark~\ref{bdd-dist} implies that $\{ k(g_n x_n) \}$ is bounded. Since $\{k(g_n)\}$ is
convergent (hence bounded), and $k$~is a homomorphism, this implies that $\{k(x_n)\}$ is
bounded. Since $k|_{X_1}$ is a homeomorphism onto $k(X_1)$, we conclude that $\{x_n\}$
is bounded. Thus, $\{ (g_n, x_n, t_n) \}$ has a convergent subsequence, as desired.

\begin{step} \label{T1-T2}
 We may assume $\delta(g) \in T_1 T_2$, for all $g \in G_1$.
 \end{step}
 Being left $\Gamma_1$-equivariant, the homeomorphism~$\phi$ factors
through to a homeomorphism $\phi^* \colon \M1 \to \M1$. By replacing $f$ with the
composition $f \circ (\phi^*)^{-1}$, we may assume $\tf(g) \in k(g) T_1 T_2$,
as desired.

\begin{defn} \label{tau-defn}
 For $x \in X_1$, because $YT_1 = k(X_1) T_1$, there is some  $\ell1(x) \in T_1$,
such that $k(x) \in Y \ell1(x)$. This means $k(x) \ell1(x)^{-1} \in Y$, so,
because $YT_2 = X_2 T_2$, there is some $\ell2(x) \in T_2$, such that $k(x)
\ell1(x)^{-1} \ell2(x) \in X_2$. Note that $\ell1(x)$ and  $\ell2(x)$ are
uniquely determined by~$x$, because $Y \cap T_1 = e$ and $X_2 \cap T_2 = e$
\see{sc-no-cpct}.
 \end{defn}

 \begin{step} \label{[Y,T]}
 $[Y,T_1]$ and $[Y,T_2]$ are contained in $Y \cap X_2$.
 \end{step}
 Because $T_1$ and $T_2$ normalize~$Y$, we know $[Y,T_1], [Y,T_2] \subset Y$, so
the proof can be completed by showing $[\Gs, \Gs] \cap Y \subset X_2$. For any $y \in
[\Gs, \Gs] \cap Y$, we know $y \in Y$, so there is some $t \in T_2$ with $yt \in X_2
\subset G_2$. However, since we also have $y \in [\Gs,\Gs] \subset G_2$
\see{ell-prop}, this implies $t = y^{-1} (yt) \in G_2$. Since $T_2 \cap G_2 = e$
\see{sc-no-cpct}, this implies $t = e$. Therefore, $y = yt \in X_2$.

\begin{step} \label{ell-equi}
 For all $g \in G_1$ and $x \in X_1$, we have $\delta(gx) = \ell1(x)^{-1}
\delta(g) \ell2(x) $.
 \end{step}
 We have
 $$ k(x) \delta(gx) = k(g)^{-1} k(gx) \delta(gx) = k(g)^{-1} \tf(gx)
 = k(g)^{-1} \tf(g) x_2 
 = \delta(g) x_2,$$
 for some $x_2 \in X_2$. Write $x_2 = yt$, with $y \in Y$ and $t \in T_2$. Then
 $$ \bigl( k(x) \ell1(x)^{-1} \bigr) \bigl( \ell1(x) \delta(gx) \bigr) 
 = k(x) \delta(gx) = \delta(g) yt
 = \bigl( \delta(g) y \delta(g)^{-1} \bigr) \delta(g) t ,$$
 so, because $\nu$ is injective, we conclude that 
 \begin{equation}
 k(x) \ell1(x)^{-1} = \delta(g) y \delta(g)^{-1}
 \label{Ypart}
 \end{equation}
 and 
 \begin{equation}
 \ell1(x) \delta(gx) = \delta(g) t .
 \label{Tpart}
 \end{equation}
 \eqref{Ypart} (and the definition of $\ell2(x)$) implies  
 $\bigl( \delta(g) y \delta(g)^{-1} \bigr) \ell2(x) \in X_2 $,
 and Step~\ref{[Y,T]} implies $[Y, T_2 T_1] \subset X_2$, 
 so we have
 $$y \ell2(x) \in [Y, T_2 T_1] \bigl( \delta(g) y \delta(g)^{-1} \bigr) \ell2(x) \subset
[Y, T_2 T_1] X_2 = X_2 ,$$
 so we must have $t = \ell2(x)$. Thus, \eqref{Tpart} yields $\ell1(x) \delta(gx) =
\delta(g) \ell2(x)$.

\begin{step} \label{delta-is-homo}
 There are homomorphisms $\delta_1 \colon G_1 \to T_1$ and $\delta_2 \colon G_1 \to
T_2$, such that $\delta(g) = \delta_1(g) \delta_2(g)$, for all $G \in G_1$, and $(\ker
\ell2)^\circ \subset \ker \delta_2$.
 \end{step}
 Let $D = \{\, (t,t) \mid t \in T_1 \cap T_2 \,\}$. Note that $D$ is a subgroup of $T_1
\times T_2$, and let $\mu \colon D \backslash (T_1 \times T_2) \to T_1 T_2$ be the
homeomorphism defined by $\mu(t_1,t_2) = t_1^{-1} t_2$. 
 Define a map $\ell{} \colon X_1 \to D \backslash (T_1 \times T_2)$ by
 $\ell{}(x) = D \cdot \bigl( \ell1(x) , \ell2(x) \bigr) $. 
 From Step~\ref{ell-equi}, we see that $\mu^{-1} \bigl( \delta(gx) \bigr) = \mu^{-1} 
\bigl( \delta(g) \bigr) \cdot \ell{}(x)$, for all $g \in G$ and $x \in X_1$. Therefore,
Lemma~\ref{torus-quot} implies that there is a homomorphism $\cover{\delta} \colon G_1
\to D \backslash (T_1 \times T_2)$, such that $\mu^{-1} \circ \delta = \cover{\delta}$. 

Let $T_2'$ be a subtorus of~$T_2$ that is complementary to $(T_1 \cap T_2)^\circ$. Then
$T_1 T_2' = T_1 T_2$ and $T_1 \cap T_2'$ is finite, so $T_1 \times T_2'$ is a finite
cover of $D \backslash (T_1 \times T_2)$. (A covering map~$\sigma$ can be defined by
$\sigma(t_1,t_2) = D \cdot (t_1^{-1},t_2)$.)  Since $G_1$ is simply connected, the
homomorphism~$\cover{\delta}$ lifts to a map into $T_1 \times T_2'$. In other words,
there are homomorphisms $\delta_1 \colon G_1 \to T_1$ and $\delta_2 \colon G_1 \to T_2'$
with $\cover{\delta}(g) = D \cdot \bigl( \delta_1(g)^{-1} , \delta_2(g) \bigr)$. Then
$\delta(g) = \mu \bigl( \cover{\delta}(g) \bigr) = \delta_1(g) \delta_2(g)$.

For $x \in \ker \ell2$, from Step~\ref{ell-equi}, we have 
 $\delta(x) = \ell1(x)^{-1} \delta(e) \ell2(x) = \ell1(x)^{-1} \in T_1$.
 Since $\delta(x) = \delta_1(x) \delta_2(x)$, and $\delta_1(x) \in T_1$, this implies
$\delta_2(x) \in T_1$. Therefore $\delta_2 ( \ker \ell2 ) \subset T_1 \cap T_2'$ is
finite, so $\delta_2 \bigl( (\ker \ell2)^\circ \bigr) = e$.

\begin{step} \label{f-is-type-C}
 $f$ is of the type described in Example~\ref{map-egs}(\ref{map-doublycrossed}).
 \end{step}
 From the definition of~$\delta$ (Step~\ref{delta-in-N_G(Y)}) and
Step~\ref{delta-is-homo}, we have
 $\tf(g) = k(g)  \cdot \delta(g) = k(g) \cdot \delta_1(g) \cdot \delta_2(g)$. 
 Furthermore, for every $\gamma \in \Gamma_1$, because $\delta(\gamma g) = \delta(g)$,
we have 
 $$\delta_1(\gamma) \cdot \delta_2(\gamma) = \delta(\gamma) = \delta(\gamma e) =
\delta(e) = e.$$
 Therefore, in order to show that $f$ arises from the construction of
Example~\ref{map-egs}\ref{map-doublycrossed}, we only need to show that the conditions
of Lemma~\ref{phi(G)isgroup} are satisfied.

Since $\addtorus{k(X_1)} = k(X_1) \prodsemi T_1$, and
$[\addtorus{k(X_1)},\addtorus{k(X_1)}] = [k(X_1),k(X_1)]$ \see{ell-prop}, we have
 \begin{align*}
 [k(X_1) , T_1]
 &\subset [\addtorus{k(X_1)},\addtorus{k(X_1)}] = [k(X_1),k(X_1)] = k([X_1,X_1]) \\
 & \subset k(X_1 \cap [G_1,G_1]) \subset k(X_1 \cap \ker \delta_1 \cap \ker \delta_2) .
 \end{align*}
 Similarly, we have $[k(G_1) , T_1] \subset k(\ker \delta_1 \cap \ker \delta_2) $.

For any $x \in X_1$, we have 
 $$\phi_1(x) \delta_2(x) = k(x) \delta_1(x) \delta_2(x) = \tf(x) \in X_2 \subset Y T_2 $$
 and $\delta_2(x) \in T_2$, so $\phi_1(X_1) \subset Y T_2$. Therefore, from
Step~\ref{[Y,T]}, we have 
 $$ [\phi_1(X_1),T_2] \subset [Y T_2, T_2] = [Y,T_2] \subset (Y \cap
X_2)^\circ .$$

Define a map $\xi \colon X_1 \to Y$ by $\xi(x) = k(x)
\ell1(x)^{-1}$. Then $\xi$ is bijective, because $k(X_1) T_1 = Y T_1$ and $k(X_1) \cap
T_1 = e$. It is also proper, because $k|_{X_1}$ is proper and $T_1$ is compact. Therefore
$\xi$ is a homeomorphism \see{proper-homeo}.

 Let $y \in (Y \cap X_2)^\circ$. Since $y \in Y$, there is some $x \in X_1$, such that
$\xi(x) = y$. Since $y \in X_2$, we have $\ell2(x) = e$, so $x \in \ker
\ell2$. Because $\xi$ is a homeomorphism, we know $\xi^{-1} \bigl( (Y \cap X_2)^\circ
\bigr)$ is connected, so we conclude that $x \in (\ker \ell2)^\circ \subset \ker
\delta_2$ (see Step~\ref{delta-is-homo}). Also, from Step~\ref{ell-equi}, we have
$\ell1(x)^{-1} \ell2(x) = \delta(x)$. Therefore 
 $$ \ell1(x)^{-1} = \ell1(x)^{-1} \cdot e = \ell1(x)^{-1} \ell2(x) = \delta(x) =
\delta_1(x) \delta_2(x) = \delta_1(x) \cdot e = \delta_1(x) ,$$
 so 
 $$ y = k(x) \ell1(x)^{-1} = k(x) \delta_1(x) = \phi_1(x) \in \phi_1(X_1) .$$
 Therefore $[\phi_1(X_1),T_2] \subset \phi_1(X_1)$, which means $T_2$ normalizes
$\phi_1(X_1)$.
 \end{proof}

\section{Non-solvable groups}

\begin{defn}[{\cite[Defn.~9.1]{Witte-super}}] \label{super-defn}
 A lattice~$\Gamma $ in a connected Lie group~$G$ is \emph{superrigid} if, for
every homomorphism $\alpha  \colon \Gamma \to \GL_n(\real )$, such that
$\closure{\Gamma ^\alpha }$ has no nontrivial, connected, compact, semisimple,
normal subgroups, there is a continuous homomorphism $\beta  \colon G \to
\closure{\Gamma ^\alpha}$, such that $\beta$ agrees with~$\alpha$ on a
finite-index subgroup of~$\Gamma $.
 \end{defn}

\begin{rem} \label{only-T}
 In the context of Definition~\ref{super-defn}, suppose $G_2$ is any connected Lie
subgroup of $\GL_n(\real)$ that contains~$\Gamma ^\alpha$. Then $\beta $ induces a
homomorphism $\quot{\beta } \colon G_2 \to \closure{G_2}/ G_2$. Since
$\closure{G_2}/ G_2$ is abelian, and $\quot{\beta }$ is trivial on a finite-index
subgroup of the lattice~$\Gamma $, we see that $\quot{\beta }(G)$ is compact and
abelian. Therefore $\beta (G) \subset G_2 S$, for any maximal compact torus~$S$ of
$\closure{\Rad G_2}$. Therefore, letting $\Gs = G_2 \prodsemi T$, for any maximal
compact torus~$T$ of $\Zi2{\Rad G_2}$, we see that there is a homomorphism $k
\colon G \to G_2 \prodsemi T$, such that $k$ agrees with~$\alpha $ on a
finite-index subgroup of~$\Gamma $.
 \end{rem}

\begin{defn}
 Let us say that a connected Lie group~$G$ is \emph{almost linear} if there is
a continuous homomorphism $\beta  \colon G \to \GL_n(\real)$, for some~$n$, such
that the kernel of~$\beta $ is finite.
 \end{defn}

The following two theorems combine to show that many lattices are superrigid.
Furthermore, by considering induced representations, it is easy to see that every
finite-index subgroup of a superrigid lattice is superrigid.

\begin{thm}[{Margulis \cite[Thm.~IX.5.12(ii), p.~327]{MargBook}}] \label{Marg-thm}
 Let $G$ be a simply connected, almost linear, semisimple Lie group, such that
$\operatorname{\real-rank}(L) \ge 2$, for every simple factor~$L$ of~$G$. Then
every lattice~$\Gamma $ in~$G$ is superrigid.
 \end{thm}

\begin{thm}[{\cite[Thm.~9.9]{Witte-super}}]
 Let $\Gamma $ be a lattice in a simply connected, almost linear Lie
group~$G$, and assume that $G$ has no nontrivial, connected, compact, semisimple,
normal subgroups. Then $\Gamma $ is superrigid iff
 \begin{itemize}
 \item $\Z{\Gamma } = \G$; and
 \item the image of~$\Gamma $ in $G/\Rad G$ is a superrigid lattice.
 \end{itemize}
 \end{thm}

\begin{thm} \label{nonsolv-fol-map}
 For $i=1,2$, let $X_i$ be a closed, unimodular subgroup of a simply connected,
almost linear Lie group~$G_i$. Assume $\Rad X_i$ is simply connected, and that
$X_i$ has no nontrivial, compact, semisimple quotients. Also assume there is no
connected, closed, normal subgroup~$N$ of $[X_1,X_1]$, such that
$\closure{\Ad_{[X_1,X_1]/N} X_1}$ is compact and nontrivial.

 For $i=1,2$, let $\Gamma_i$ be a lattice in $G_i$. Assume that $G_i$ has no
nontrivial, connected, compact, semisimple, normal subgroups, that $\Zi{i}{\Gamma
_i} = \Gi{i}$, and that $\Gamma _1$ is superrigid in~$G_1$.

Assume, furthermore, that the foliation of~$\Mp1$ by cosets of~$X_1$
has a dense leaf, for every finite-index subgroup $\Gamma_1'$ of~$\Gamma_1$.

 Let~$f \colon \Mp1 \to \M2$ be a homeomorphism, such that $f$ maps each leaf of
the foliation of~$\M1$ by cosets of~$X_1$ onto a leaf of the foliation of~$\M2$
by cosets of~$X_2$.

 Then, for some finite-index subgroup $\Gamma_1'$ of~$\Gamma_1$, there exists 
 \begin{itemize}
    \item a map $a \colon \Mp1 \to \Mp1$ of type~\ref{map-egs}\ref{map-leaftoself};
    \item a map $b \colon \M2 \to \M2$ of type~\ref{map-egs}\ref{map-translate};
and
    \item a map $c \colon \Mp1 \to \M2$ of
type~\ref{map-egs}\ref{map-doublycrossed},
  \end{itemize}
 such that $f(\Gamma _1 g) = b\Bigl( c \bigl( a(\Gamma_1' g) \bigr) \Bigr)$,
 for all $g \in G_1$.

In the definition of~$c$, we may take $\Gs$ as defined in
Remark~\ref{only-T}. We may take $T_1$ to be the elliptic part of~$k(X_1)$
and let $T_2$ be an appropriately chosen elliptic part of~$r^{-1} X_2 r$, where
$k$ is the homomorphism used in the construction of~$c$, and $r$ is the element
of~$G_2$ used in the construction of~$b$.
 \end{thm}

\begin{proof}[Sketch of proof]
 The proof of Theorem~\ref{solv-fol-map} applies with only minor changes; we point
out the substantial differences.

A change is required already in the first paragraph of the proof. Assume
$f(\Gamma _1) = \Gamma _2$. Then $f$ lifts to a homeomorphism $\tf \colon G_1
\to G_2$ with $\tf(e) = e$. Since $G_2$ is almost linear, $\Gamma _1$ is
superrigid, and $G_1$ is simply connected, it is not difficult to see from
Remark~\ref{only-T} that there is a finite-index subgroup~$\Gamma _1'$ of~$\Gamma
_1$, such that $\tf|_{\Gamma _1'}$ extends to a homomorphism $k \colon G_1 \to
\Gs$. For simplicity, replace $\Gamma _1$ with~$\Gamma _1'$.

See \cite[Defns.~4.3 and~4.8]{Witte-fols} for the definition of elliptic and
nonelliptic parts of a subgroup of~$\Gs$. Note that the assumptions on~$X_i$ imply
that $T_i$ is a torus, and $X_i \cap T_i$ is finite. Also, $X_i$ has no nontrivial
connected, compact, normal subgroups.

Because $\M1$ may not be compact, the second sentence of Remark~\ref{bdd-dist} may
not be valid, so the arguments of Steps~\ref{same-nonell}
and~\ref{delta-in-N_G(Y)} need to be modified, as in Steps~3 and~4 of the proof of 
\cite[Thm.~6.1]{Witte-fols}.

The conclusion of Step~\ref{k-homeo} should be weakened slightly: instead of being
a homeomorphism, the restriction of~$k$ to~$X_1$ is a finite-to-one covering map.
Similarly, the maps $\nu$, $\nu_1$, and~$\nu_2$ of Step~\ref{nu-homeo} are
finite-to-one covering maps. For example, to see this in the case of~$\nu$, let
$T_1'$ be a subtorus of~$T_1$, such that $T_1' \cap T_2$ is finite, and $T_1 = T_1'
(T_1 \cap T_2)$. Then the group $(Y \prodsemi T_1') \times T_2$ acts on $G_2$ by 
 $(y,t_1,t_2) \cdot g = yt_1gt_2$.
 The orbit of~$e$ under this action is $Y T_1 T_2$, and the stabilizer of~$e$ is
finite. So the map $(y,t_1,t_2) \mapsto y t_1 t_2$ is a covering map with finite
fibers. The space $Y \times T_1 T_2 \homeo \bigl( (Y \prodsemi T_1') \times T_2
\bigr) / (T_1' \cap T_2)$ is an intermediate covering space, with covering
map~$\nu$.

In the proof of Step~\ref{phi-defn}, $\nu_1^{-1} \bigl( \delta(g) \bigr)$ may not
be a single point, but, because $G_1$ is simply connected, there is a continuous
function $\hat{\chi} \colon G_1 \to X_1  \times T_1 T_2$ with $\nu_1 \bigl(
\hat{\chi}(g) \bigr) = \delta(g)$ (and $\hat{\chi}(e) = e$). Define $\chi(g)$ to
be the first component of $\hat{\chi}(g)$. (A similar device is used to
define~$\chi'$.) Because $\nu_1$ is finite-to-one, and the lift $\hat{\chi}$ is
determined by its value at any one point, there is a finite-index subgroup~$\Gamma$
of~$\Gamma _1$, such that $\chi(\gamma g) = \chi(g)$, for all $g \in G_1$ and
$\gamma \in \Gamma $. Then, by replacing~$\Gamma _1$ with~$\Gamma $, we may assume
$\phi$ is left $\Gamma _1$-equivariant.

Still in the proof of Step~\ref{phi-defn}, to see that $\zeta$ is a
finite-to-one covering map, note that the map $\zeta_1 \colon G_1 \times X_1 \to
G_1 \times X_2$ defined by $\zeta_1(g,x) = \bigl( g, \tf(g)^{-1} \tf(gx) \bigr)$
is a homeomorphism (by the argument in the last paragraph of Step~\ref{phi-defn}),
and we have $\zeta(g,x,t) = (\operatorname{Id} \times \nu_2) \bigl(\zeta_1(g,x),t
\bigr)$. 

Also in the proof of Step~\ref{nu-homeo}, let us show that $\psi$ is the inverse
of~$\phi$. For all $g \in G_1$, we have 
 $$\tf(g) \in k \bigl( \phi(g) \bigr) T_1 T_2
 \qquad \text{and} \qquad
 \tf \Bigl( \psi \bigl( \phi(g) \bigr) \Bigr) \in k \bigl( \phi(g) \bigr) T_1
T_2 .$$
 Therefore, the two maps $\lambda _1,\lambda _2 \colon G_1 \to X_2 \times T_2 T_1$
defined by
 $$ \lambda _1(g) = \Bigl( e, \tf(g)^{-1} k \bigl( \phi(g) \bigr) \Bigr) $$
  and
 $$ \lambda _2(g) = \biggl( \tf(g)^{-1}\tf \Bigl( \psi \bigl( \phi(g) \bigr) \Bigr)
, \tf \Bigl( \psi \bigl( \phi(g) \bigr) \Bigr)^{-1} k \bigl( \phi(g) \bigr)
\biggr) $$
 both satisfy $\nu_2 \bigl( \lambda _i(g) \bigr) = \tf(g)^{-1} k \bigl( \phi(g)
\bigr)$. Since $\lambda _1(e) = (e,e) = \lambda _2(e)$, and $\nu_2$ is a covering
map, this implies $\lambda _1 = \lambda _2$. By comparing the first coordinates,
and noting that $\tf$ is a homeomorphism, we conclude that $\psi \circ \phi$ is
the identity map. Similarly, $\phi \circ \psi$ is also the identity.

In Definition~\ref{tau-defn}, since $X_i \cap T_i$ may be nontrivial, the
functions $\ell1$ and~$\ell2$ may not be well defined. However, there is a finite
cover~$\cover{X_1}$ of~$X_1$, such that $\ell{i}$ is a well-defined map
from~$\cover{X_1}$ to~$T_i$ with $\ell{i}(e) = e$.

The conclusion of Step~\ref{[Y,T]} can be established as follows.  Because $T_2$
normalizes~$Y$, we have $[Y,T_2] \subset Y$. Furthermore, because $T_2$ is abelian and
is in $\Zi2{X_2}$, we have $[T_2 X_2, T_2 X_2] = [X_2, X_2] \subset X_2$. Therefore,
$[Y,T_2] \subset Y \cap X_2$.  Since $[X_1,X_1] \subset Y$, we have
$\Ad_{[X_1,X_1]/[Y,Y]} Y = e$, so $\closure{\Ad_{[X_1,X_1]/[Y,Y]} X_1} =
\Ad_{[X_1,X_1]/[Y,Y]} T_1$ is compact. Therefore, by hypothesis, it must be trivial.
Since $\Ad_{T_1 X_1/[X_1,X_1]} T_1$ is also trivial, and $\Ad_{X_1} T_1$ is semisimple,
this implies $[T_1 X_1, T_1] \subset [Y,Y] \subset X_2$. Therefore, $[Y, T_1] \subset
Y \cap X_2$.

The argument of Step~\ref{ell-equi} shows that, for all $g \in G_1$ and $x \in
\cover{X_1}$, we have
 $$ \delta(gx) \in (T_1 \cap Y) \ell1(x)^{-1} \delta(g) \ell2(x) (T_2 \cap Y) .$$
 (This calculation uses the observation that $T_1 \cap Y$ and $T_2 \cap Y$ are
contained in $Y \cap X_2$. To see this, note that $T_i \cap Y$ is contained in a
maximal compact subgroup~$K$ of~$Y$. Because $\Rad Y$ is simply connected, we see
that $K$ is contained in a Levi subgroup of~$Y$, so $K \subset [Y,Y] \subset Y
\cap X_2$.) Since $\delta$, $\ell1$, and~$\ell2$ are continuous, and $T_i \cap Y$
is finite, this implies the equation in the conclusion of Step~\ref{ell-equi}.

For Step~\ref{delta-is-homo}, it is necessary to prove a slightly modified version of
Lemma~\ref{torus-quot}.

In Step~\ref{f-is-type-C}, although the map $\xi \colon \cover{X_1} \to Y$ is not a
homeomorphism, we must have $(Y \cap X_2)^\circ \subset \xi\bigl( (\ker \ell2)^\circ
\bigr)$, because $\xi$ is a covering map.
 \end{proof}

\end{document}